\magnification=1200
\overfullrule=0pt
\centerline {\bf A strict minimax inequality criterion and some of its
consequences}\par
\bigskip
\bigskip
\centerline {BIAGIO RICCERI}\par
\bigskip
\bigskip
\noindent
{\bf Abstract:} In this paper, we point out a very flexible scheme
within which a strict minimax inequality occurs. We then show the
fruitfulness of this approach presenting a series of various consequences.
Here is one of them: \par
Let $Y$ be a
finite-dimensional real Hilbert space,
$J:Y\to {\bf R}$ a $C^1$ function with locally Lipschitzian derivative,
and $\varphi:Y\to [0,+\infty[$ a $C^1$ convex function with locally
Lipschitzian derivative at $0$ and $\varphi^{-1}(0)=\{0\}$.\par
Then, for each $x_0\in Y$ for wich $J'(x_0)\neq 0$, there exists $\delta>0$
such that, for each $r\in ]0,\delta[$, the restriction of $J$ to
$B(x_0,r)$ has a unique global minimum $u_r$ which
satisfies
$$J(u_r)\leq J(x)-\varphi(x-u_r)$$
for all $x\in B(x_0,r)$,
where $B(x_0,r)=\{x\in Y: \|x-x_0\|\leq r\}\ .$ \par
\par
\bigskip
\bigskip
\bigskip
{\bf 1. Introduction}\par
\bigskip
Let $X, \Lambda$ be two non-empty sets and $f:X\times \Lambda\to {\bf R}$ a given function.
Clearly, we have
$$\sup_{\lambda\in \Lambda}\inf_{x\in X}f(x,\lambda)\leq 
\inf_{x\in X}\sup_{\lambda\in \Lambda}f(x,\lambda)\ .$$
The classical minimax theory deals with minimax theorems, that is to say with results
which ensure the validity of the equality
$$\sup_{\lambda\in \Lambda}\inf_{x\in X}f(x,\lambda)=
\inf_{x\in X}\sup_{\lambda\in \Lambda}f(x,\lambda)\ .$$
In this paper, to the contrary, we are interested in the strict inequality
$$\sup_{\lambda\in \Lambda}\inf_{x\in X}f(x,\lambda)<
\inf_{x\in X}\sup_{\lambda\in \Lambda}f(x,\lambda)\ .\eqno{(1})$$
While the importance of minimax theorems in various fields
is well established
since long time, the interest in the study of the strict minimax
inequality is much
more recent, being been originated in
[12] and [5]. Actually, in [12], the following result was obtained:\par
\medskip
THEOREM A. - {\it  Let $X$ be a compact metric space, $\Lambda\subset
{\bf R}$ a compact interval and $f:X\times \Lambda\to {\bf R}$ a
function which is lower semicontinuous in $X$ and continuous and concave
in $\Lambda$. Assume also that $(1)$ holds.\par
Then, there exists an open interval $A\subset \Lambda$ such that, for
each $\lambda\in A$, the function $f(\cdot,\lambda)$ has a local, not
global minimum.}\par
\medskip
In turn, in [5], Theorem A was used to get the following\par
\medskip
THEOREM B. - {\it Let $X$ be a separable and reflexive real
Banach space, $\Lambda\subseteq {\bf R}$ an interval, and
$f:X\times \Lambda\to {\bf R}$ a continuous
function which is concave in $\Lambda$, and
sequentially weakly lower semicontinuous,
coercive, and satisfying the Palais-Smale condition in $X$.
Assume also that $(1)$ holds.\par
Then, there exist an open interval $A\subseteq \Lambda$
and a positive real number $\rho$, such that, for each $\lambda
\in \Lambda$, the equation
$$f'_{x}(x,\lambda)=0$$
has at least three solutions in $X$ whose norms are less than $\rho$.}\par
\medskip
When $f(x,\cdot)$ is affine, Theorem B 
is one of the most often applied multiplicty results in dealing with
problems of variational nature.
In this connection, we refer to [8] for a comprehensive account of the relevant
literature.\par
\smallskip
A further very recent motivation for the study of $(1)$ comes from [9] and
[3]. In particular, in [3], S. J. N. Mosconi proved the following\par
\medskip
THEOREM C. - {\it Let $X$ be a topological space, $\Lambda\subseteq
{\bf R}^n$ a convex set with non-empty interior, $f:X\times \Lambda\to
{\bf R}$ a function which is lower semicontinuous and inf-compact in $X$,
and upper semicontinuous and concave in $\Lambda$. Assume also that
$(1)$ holds.\par
Then, there exists $\lambda^*\in \hbox {\rm int}(\Lambda)$ such that
the function $f(\cdot,\lambda^*)$ has at least two global minima.}\par
\medskip
More precisely, Theorem C extends to the case $n>1$ a result previously
established in [9] (for $n=1$) under
the less restrictive assumption of quasi-concavity for $f(x,\cdot)$
(see also [14]). In turn, in Section 3, we will extend Theorem C to the
case where $\Lambda$ is a convex set in an arbitrary Hausdorff topological
vector space (Theorem 3.2).\par
\smallskip
It is worth noticing that Theorem A, to the contrary, does not admit such
an extension. Indeed, take $X=\{(t,s)\in {\bf R}^2: t^2+s^2=1\}$,
$\Lambda=\{(\mu,\nu)\in {\bf R}^2 : \mu^2+\nu^2\leq 1\}$ and
$$f(t,s,\mu,\nu)=\mu t+\nu s$$
for all $(t,s,\mu,\nu)\in X\times\Lambda$. Of course, we have
$$\sup_{(\mu,\nu)\in \Lambda}\inf_{(t,s)\in X}f(t,s,\mu,\nu)=0<1=
\inf_{(t,s)\in X}\sup_{(\mu,\nu)\in \Lambda}f(t,s,\mu,\nu)\ .$$
Moreover, the function $f$ is continuous in $X\times\Lambda$ and
concave in $\Lambda$. Nevertheless, for each $(\mu,\nu)\in \Lambda$,
the function $f(\cdot,\cdot,\mu,\nu)$ has no local, not global minima in
$X$.\par
\smallskip
The aim of this paper is to establish Theorem 3.1 below (whose formulation
has been inspired by Proposition 1 of [10]) which provides a very
flexible scheme within which $(1)$ is obtained.\smallskip
\smallskip
Just to stress the great flexibility of Theorem 3.1, we now state six
of its consequences, of a quite different nature, whose proofs will be given in
Section 3.\par
\smallskip
First, we fix a few relevant definitions.\par
\smallskip
A non-empty set $C$ in a normed space $S$ is said to be
uniquely remotal with respect to a set $D\subseteq S$ if, for each $y\in D$, there exists a unique $x\in C$
such that
$$\|x-y\|=\sup_{u\in C}\|u-y\|\ .$$
The main problem in theory of such sets is to know if they are singletons.\par
\smallskip
If $E, F$ are two real vector spaces and $D$ is a convex subset of $E$, 
we say that an operator $\Phi:D\to F$ is affine if
$$\Phi(\lambda x+(1-\lambda)y)=\lambda \Phi(x)+(1-\lambda)\Phi(y)$$
for all $x, y\in D$, $\lambda\in [0,1]$.
\smallskip
Let $(T,{\cal F},\mu)$ be a measure space, 
 $E$  a real Banach space and $p\geq 1$.
\par
\smallskip
As usual, $L^{p}(T,E)$  denotes the space of all (equivalence
classes of) strongly $\mu$-measurable functions $u : T\rightarrow E$ 
such that
$\int_{T}\parallel u(t)\parallel^{p} d\mu<+\infty$, equipped with
the norm $$\parallel u\parallel_{L^{p}(T,E)}=
\left ( \int_{T}\parallel u(t)\parallel^{p}d\mu\right ) ^{1\over p}\ .$$
\smallskip
A set $D\subseteq L^{p}(T,E)$ is said to be decomposable if, for
every $u,v\in D$ and every $A\in {\cal F}$, the function
 $$t\to
\chi_{A}(t)u(t)+(1-\chi_{A}(t))v(t)$$ belongs to $D$, where $\chi_{A}$
denotes the characteristic function of $A$.\par
\smallskip
Here is the series of statements announced above.\par
\medskip
THEOREM 1.1. - {\it Let $Y$ be a real normed space and let $X\subseteq Y$
be a non-empty compact uniquely
remotal set with respect to \hbox {\rm conv}$(X)$.\par
Then, $X$ is a singleton.}\par
\medskip
THEOREM 1.2. - {\it Let $X$ be a finite-dimensional real Hilbert space and 
$J:X\to {\bf R}$ a $C^1$ function. Set
$$\eta=\liminf_{\|x\|\to +\infty}{{J(x)}\over {\|x\|^2}}$$
and
$$\theta=\inf\left \{
{{J(x)-J(u)}\over {\|x-u\|^2}} : (u,x)\in M_{J}\times
X\hskip 3pt \hbox {\rm with}\hskip 3pt x\neq u\right
\}$$
where $M_J$ denotes the set of all global minima of $J$.
Assume that $$\theta<\eta\ .$$
Then, for each $\mu\in ]2\theta,2\eta[$, there exists $y_{\mu}\in X$ such that
the equation
$$J'(x)-\mu x=y_{\mu}$$
has at least three solutions.}\par
\medskip
THEOREM 1.3. - {\it Let $Y$ be a 
finite-dimensional real Hilbert space,
$J:Y\to {\bf R}$ a $C^1$ function with locally Lipschitzian derivative,
and $\varphi:Y\to [0,+\infty[$ a $C^1$ convex function with locally
lipschitzian derivative at $0$ and $\varphi^{-1}(0)=\{0\}$.\par
Then, for each $x_0\in Y$ for wich $J'(x_0)\neq 0$, there exists $\delta>0$
such that, for each $r\in ]0,\delta[$, the restriction of $J$ to
$B(x_0,r)$ has a unique global minimum $u_r$ which
satisfies
$$J(u_r)\leq J(x)-\varphi(x-u_r)$$
for all $x\in B(x_0,r)$,
where
$$B(x_0,r)=\{x\in Y: \|x-x_0\|\leq r\}\ .$$}
\medskip
THEOREM 1.4. - {\it Let $X$ be a convex subset of a real vector space, 
$(Y,\langle\cdot,\cdot\rangle)$ a real inner product space,
 $I:X\to {\bf R}$ a bounded above convex function and $\Phi:X\to Y$ an affine
 operator such that
$$\inf_{z\in \Phi(X)}\langle z,y\rangle=-\infty$$
for all $y\in \Phi(X)\setminus \{0\}$.\par
Then, the set of all global minima of the function
$x\to I(x)+\|\Phi(x)\|^2$ is contained in
 $\Phi^{-1}(0)$.}
\medskip
THEOREM 1.5. - {\it Let $(T,{\cal F},\mu)$ be a non-atomic
measure space, with $0<\mu(T)<+\infty$, $E$ a real Banach space, $f, g:E\to {\bf R}$ two
functions, with $f$ lower
semicontinuous and $g$ continuous, satisfying
$$\sup_Ef<+\infty\ , \inf_Eg=-\infty\ , \sup_Eg=+\infty$$
and
$$\max\left \{ \sup_{x\in E}{{|f(x)|}\over {1+\|x\|^p}},
\sup_{x\in E}{{|g(x)|}\over {1+\|x\|^p}}
\right \} <+\infty$$
for some $p\geq 1$.
Finally, assume that $f$ has at most one global minimum and that no zero
of $g$ is a global minimum of $f$. For each $u\in L^p(T,E)$, set
$$J(u)=\int_Tf(u(t))d\mu+\left ( \int_Tg(u(t))d\mu\right )^2\ .$$
Then, the restriction of
the functional $J$
to any decomposable subset of $L^p(T,E)$ containing the constant
functions has no global
minima.}\par
\medskip
THEOREM 1.6. - {\it Let $X$ be a convex subset of a real vector space,
$I:X\to {\bf R}$ a convex function and $P:X\to ]0,+\infty[$ a concave 
 function.\par
Then, for each $\mu>0$, a point $u\in X$ is a global minimum of the
function $x\to I(x)-\mu\log P(x)$ if and only if one has
$$I(u)\leq I(x)-\mu\left ( {{P(x)}\over {P(u)}}-1\right )$$
for all $x\in X$.}\par
\bigskip
{\bf 2. Notations}\par
\bigskip
To state our results in a more compact form, we now fix some notations.\par
\smallskip
Here and in the sequel, $X$ is a non-empty set, $\Lambda, Y$ are two topological
spaces, $y_0$ is a point in $Y$.\par
\smallskip
A family ${\cal N}$ of non-empty subsets of $X$ is said to be a filtering
cover of $X$ if $\cup_{A\in {\cal N}}A=X$ and for each $A_1, A_2\in
{\cal N}$ there is $A_3\in {\cal N}$ such that $A_1\cup A_2\subseteq
A_3$.\par
\smallskip
We denote by ${\cal G}$
the family of all lower semicontinuous functions $\varphi:Y\to [0,+\infty[$,
with $\varphi^{-1}(0)=\{y_0\}$, such that, for each neighbourhood $V$ of $y_0$,
one has
$$\inf_{Y\setminus V}\varphi>0\ .\eqno{(2)}$$
Moreover, we denote by ${\cal H}$ the family of all functions
$\Psi:X\times\Lambda\to Y$ such that, for each $x\in X$, 
$\Psi(x,\cdot)$ is continuous, injective, open, 
takes the value $y_0$ at a point $\lambda_x$ and the function
$x\to \lambda_x$ is not constant. 
Furthermore, we denote by ${\cal M}$ the family of all
functions $J:X\to {\bf R}$ whose set of all global minima
(noted by $M_{J}$) is non-empty.\par
\smallskip
Finally, for each $\varphi\in {\cal G}$, $\Psi\in {\cal H}$ 
 and $J\in {\cal M}$, we put
$$\theta(\varphi,\Psi,J)=\inf\left \{
{{J(x)-J(u)}\over {\varphi(\Psi(x,\lambda_u))}} : (u,x)\in M_{J}\times
X\hskip 3pt \hbox {\rm with}\hskip 3pt \lambda_x\neq \lambda_u\right
\}\ .$$
\bigskip
{\bf 3. Results and proofs}\par
\bigskip
With such notations, our abstract criterion reads as follows:\par
\medskip
THEOREM 3.1. - {\it Let $\varphi\in {\cal G}$, $\Psi\in {\cal H}$  
 and $J\in {\cal M}$. \par
Then, for each $\mu>\theta(\varphi,\Psi,J)$ and each filtering cover
${\cal N}$ of $X$, there exists $A\in {\cal N}$ such that
$$\sup_{\lambda\in \Lambda}\inf_{x\in A}
(J(x)-\mu\varphi(\Psi(x,\lambda)))<
\inf_{x\in A}\sup_{z\in A}
(J(x)-\mu\varphi(\Psi(x,\lambda_z)))\ .$$}\par
\smallskip
PROOF. 
 Let $\mu>\theta(\varphi,\Psi,J)$ and let ${\cal N}$ be
a filtering cover of $X$. 
 Choose $u\in M_{J}$ and $x_1\in X$,
with $\lambda_{x_1}\neq \lambda_u$, such that
$$J(x_1)-\mu\varphi(\Psi(x_1,\lambda_u))<J(u)\ .$$
Due to the nature of ${\cal N}$, there exists $A\in
{\cal N}$ such that $u, x_{1}\in A$. We have
$$0\leq\inf_{z\in A}\varphi(\Psi(x,\lambda_z))\leq
\varphi(\Psi(x,\lambda_x))=0$$
for all $x\in A$, and so, since $u$ is a global
minimum of $J$, it follows that
$$\inf_{x\in A}\sup_{z\in A}
(J(x)-\mu\varphi(\Psi(x,\lambda_z)))=\inf_{x\in A}
\left ( J(x)-\mu\inf_{z\in A}
\varphi(\Psi(x,\lambda_z))\right )$$
$$ =\inf_XJ=J(u)\ .\eqno{(3)}$$
Since the function
$\varphi(\Psi(x_1,\cdot))$ is lower semicontinuous at
$\lambda_u$, 
there are $\epsilon>0$ and
a neighbourhood $U$ of $\lambda_u$ such that
$$J(x_1)-\mu\varphi(\Psi(x_1,\lambda))<J(u)-\epsilon$$
for all $\lambda\in U$. 
So, we have
$$\sup_{\lambda\in U}\inf_{x\in A}
(J(x)-\mu\varphi(\Psi(x,\lambda)))\leq\sup_{\lambda\in U} 
 (J(x_1)-\mu\varphi(\Psi(x_1,\lambda)))\leq J(u)-\epsilon
\ .\eqno{(4)}$$
Since $\Psi(u,\cdot)$ is open,
 the set $\Psi(u,U)$ is a neighbourhood of $y_0$.
Hence, by $(2)$, we have
$$\nu:=\inf_{y\in Y\setminus \Psi(u,U)}\varphi(y)>0\ .\eqno{(5)}$$
Moreover, since $\Psi(u,\cdot)$ is injective,
 if $\lambda\not\in U$ then $\Psi(u,\lambda)\not
\in \Psi(u,U)$. So,
 from $(4)$ and $(5)$, it follows that
$$\sup_{\lambda\in \Lambda\setminus U}\inf_{x\in A}
(J(x)-\mu\varphi(\Psi(x,\lambda))\leq J(u)-
\mu\inf_{\lambda\in \Lambda\setminus U}\varphi(\Psi(u,\lambda))\leq
J(u)-\mu\nu\ . \eqno{(6)}$$
Now, the conclusion comes directly from $(3)$, $(4)$, $(5)$ and
$(6)$.\hfill $\bigtriangleup$\par
\medskip
REMARK 3.1. - From the conclusion of Theorem 3.1 it clearly follows that,
for any set $D\subseteq\Lambda$ with $\lambda_x\in D$ for all $x\in A$,
one has
$$\sup_{\lambda\in D}\inf_{x\in A}
(J(x)-\mu\varphi(\Psi(x,\lambda)))<
\inf_{x\in A}\sup_{\lambda\in D}
(J(x)-\mu\varphi(\Psi(x,\lambda)))\ .$$
\medskip
REMARK 3.2. - From the definition of $\theta(\varphi,\Psi,J)$, it clearly follows that
$u\in M_J$ if and only if $u$ is a global minimum of the
function $x\to J(x)-\theta(\varphi,\Psi,J)\varphi(\Psi(x,\lambda_u))$.
So, when
$\theta(\varphi,\Psi,J)>0$, from knowing that
$$J(u)\leq J(x)$$
for all $x\in X$, we automatically get
$$J(u)\leq J(x)-\theta(\varphi,\Psi,J)\varphi(\Psi(x,\lambda_u))$$
for all $x\in X$, which is a much better inequality since
$\varphi(y)>0$ for all $y\in Y\setminus \{y_0\}$.\par 
\medskip
REMARK 3.3 - It is likewise important to observe that if
$\theta(\varphi,\Psi,J)>0$,
then the function $x\to \lambda_x$ is constant in $M_{J}$. As
a consequence, if $\theta(\varphi,\Psi,J)>0$ and
the function $x\to \lambda_x$ is injective,
then $J$ has a unique global minimum. In
particular, note that $x\to \lambda_x$ is injective
 when $\Psi(\cdot,\lambda)$ is injective for
all $\lambda\in \Lambda$.\par
\medskip
REMARK 3.4. - Remarks 3.2 and 3.3 show the interest in knowing
when $\theta(\varphi,\Psi,J)>0$. Theorem 3.1 can also be useful
for this. Indeed, if for some $\mu>0$, there is a filtering cover
${\cal N}$ of $X$ such that 
$$\sup_{\lambda\in \Lambda}\inf_{x\in A}
(J(x)-\mu\varphi(\Psi(x,\lambda)))\geq
\inf_{x\in A}\sup_{z\in A}
(J(x)-\mu\varphi(\Psi(x,\lambda_z)))$$
for all $A\in {\cal N}$, then $\theta(\varphi,\Psi,J)\geq \mu$.\par
\medskip
As we said in the Introduction,
we now extend Theorem C to the case where
$\Lambda$ is a convex subset of an arbitrary Hausdorff
topological vector space. Recall that, when $U$ is a topological
space, a function $\psi:U\to {\bf R}$ is said to be inf-compact
if, for each $r\in {\bf R}$, the set $\{x\in U : \psi(x)\leq r\}$ is compact.\par
\medskip
THEOREM 3.2. - {\it Let $X$ be a topological space, $E$ a real Hausdorff
topological vector space, $\Lambda\subseteq
E$ a convex set, $f:X\times \Lambda\to
{\bf R}$ a function which is lower semicontinuous and inf-compact in $X$,
and upper semicontinuous and concave in $\Lambda$. Assume also that
$(1)$ holds.\par
Then, there exists $\lambda^*\in \Lambda$ such that
the function $f(\cdot,\lambda^*)$ has at least two global minima.}\par
\medskip
PROOF. Denote by ${\cal S}_{\Lambda}$ the family of all finite-dimensional
convex subsets of $\Lambda$.
Arguing by contradiction, assume that, for each $\lambda\in \Lambda$,
the function $f(\cdot,\lambda)$ has a unique global minimum. 
Fix $S\in {\cal S}_{\Lambda}$. We claim that
$$\sup_S\inf_Xf=\inf_X\sup_Sf\ .\eqno{(7)}$$
Let $n$ be the dimension of $S$. Since $E$ is Hausdorff, there exists
an affine homeomorphism $\Phi$ between aff$(S)$, the affine hull
of $S$, and ${\bf R}^n$.
Thus, $Y:=\Phi(S)$ is a convex subset of ${\bf R}^n$ of dimension
$n$, and so its interior is non-empty. Now, put
$$\tilde f(x,y)=f(x,\Phi^{-1}(y))$$
for all $(x,y)\in X\times Y$. Of course, the function $\tilde f$ is
lower semicontinuous and inf-compact in $X$, while is upper semicontinuous
and concave in $Y$. Therefore, since $\tilde f(\cdot,y)$ has a unique
global minimum for all $y\in Y$, Theorem C ensures that
$\tilde f$ does not satisfy $(1)$, that is 
$$\sup_Y\inf_X\tilde f=\inf_X\sup_Y\tilde f\ .\eqno{(8)}$$
Of course, one has
$$\sup_Y\inf_X\tilde f=\sup_S\inf_Xf$$ 
and
$$\inf_X\sup_Y\tilde f=\inf_X\sup_Sf\ .$$
So, $(7)$ comes directly from $(8)$.
Now, since $f$ satisfies $(1)$, we can fix $r$ so that
$$\sup_{\Lambda}\inf_Xf<r<\inf_X\sup_{\Lambda}f\ .\eqno{(9)}$$
For each $S\in {\cal S}_{\Lambda}$, put
$$C_S=\bigcap_{\lambda\in S}\{x\in X : f(x,\lambda)\leq r\}\ .$$
Clearly, $C_S$ is closed. Note also that $C_S$ is non-empty.
Indeed, otherwise, we would have
$$r\leq\inf_X\sup_Sf$$
and so, by $(7)$,
$$r\leq\sup_S\inf_Xf$$
against $(9)$. Of course, if $S_1,...,S_k\in {\cal S}_{\Lambda}$,
then there is $\hat S\in {\cal S}_{\Lambda}$ such that
$$\bigcup_{i=1}^k S_i\subseteq\hat S\ .$$
So
$$C_{\hat S}\subseteq \bigcap_{i=1}^k C_{S_i}\ .$$
In other words, the family of closed compact sets $\{C_S\}_{S\in {\cal S}_{\Lambda}}$
has the finite intersection property, and so its intersection is non-empty.
Hence, there is $\hat x\in X$ such that
$$f(\hat x,\lambda)\leq r$$
for all $\lambda\in S$ and all $S\in {\cal S}_{\Lambda}$. Consequently
$$\inf_X\sup_{\Lambda}f\leq r\ ,$$
against $(9)$. Such a contradiction completes the proof.
\hfill $\bigtriangleup$\par
\medskip
A joint, direct application of Theorems 3.1 and 3.2 gives\par
\medskip
THEOREM 3.3. - {\it Let $\varphi\in {\cal G}$, $\Psi\in {\cal H}$  
 and $J\in {\cal M}$. 
Moreover, assume that $X$ is a topological space, that
$\Lambda$ is a real Hausdorff topological vector space and that
$\varphi(\Psi(x,\cdot))$ is convex for each $x\in X$. Finally, let
$\mu>\theta(\varphi,\Psi,J)$ and ${\cal N}$ be a filtering cover of
$X$
such that, for each $A\in {\cal N}$,
 the function $x\to J(x)-\mu\varphi(\Psi(x,\lambda))$ is
lower semicontinuous and inf-compact in $A$ for all $\lambda\in 
\hbox {\rm conv}(\{\lambda_x : x\in X\})$.
\par
Under such hypotheses,
there exist $A\in {\cal N}$ and $\lambda^*\in
\hbox {\rm conv}(\{\lambda_x : x\in A\})$
 such
that the restriction
of the function $x\to J(x)-\mu\varphi(\Psi(x,\lambda^*))$ to $A$
has at least two global minima.}\par
\smallskip
PROOF. For each $(x,\lambda)\in X\times {\Lambda}$, put
$$f(x,\lambda)=J(x)-\mu\varphi(\Psi(x,\lambda))\ .$$
Recalling that $\Psi(x,\cdot)$ is continuous and $\varphi$ is
lower semicontinuous, we have that, for each $x\in X$, the function
$f(x,\cdot)$ is upper semicontinuous and concave in $\Lambda$. 
By Theorem 3.1, there exists $A\in {\cal N}$ such that
$$\sup_{\lambda\in D}\inf_{x\in A}f(x,\lambda)<
\inf_{x\in A}\sup_{\lambda\in D}f(x,\lambda)\ ,$$
where
$$D=\hbox {\rm conv}(\{\lambda_x : x\in A\})\ .$$
Now, the conclusion comes directly applying Theorem 3.2 to
the restriction of $f$ to $A\times D$.\hfill $\bigtriangleup$\par
\medskip
On the basis of Theorem 3.3, we can give the proofs of Theorems 1.1-1.3.\par
\medskip
{\it Proof of Theorem 1.1}. Arguing by contradiction, assume that $X$ contains at least 
two points.
Now, apply Theorem 3.3 taking: $\Lambda=Y$, $y_0=0$, $\varphi(x)=\|x\|$, $\Psi(x,\lambda)=
x-\lambda$, $J=0$ and ${\cal N}=\{X\}$. Note that we are allowed to apply Theorem 3.3 since
$x\to \lambda_x$ is not constant. Then, it would exists
 $\lambda^*\in \hbox {\rm conv}(X)$
such that the function $x\to -\|x-\lambda^*\|$ has at least two global minima in
$X$, against the hypotheses.\hfill $\bigtriangleup$\par
\medskip
REMARK 3.5. - Observe that Theorem 1.1 improves a classical result by
V. L. Klee ([1]) under
two aspects: $Y$ does not need to be complete and $\overline {\hbox {\rm conv}}(X)$
is replaced by $\hbox {\rm conv}(X)$. Note also that our proof is completely different
from that of Klee which is based on the Schauder fixed point theorem.
\medskip
{\it Proof of Theorem 1.2}. 
Let $\mu\in ]2\theta,2\eta[$. We clearly
have
$$\lim_{\|x\|\to +\infty}\left ( J(x)-{{\mu}\over {2}}\|x-
\lambda\|^2\right ) =+\infty
\eqno{(10)}$$
for all $\lambda\in X$. So, since $X$ is finite-dimensional,
the function $x\to
J(x)-{{\mu}\over {2}}\|x-\lambda\|^2$ is continuous and inf-compact for
all $\lambda\in X$. Therefore, we can apply Theorem 3.3 
 taking: $X=Y=\Lambda$, $y_0=0$, $\varphi(y)=\|y\|^2$,
$\Psi(x,\lambda)=x-\lambda$ and ${\cal N}=\{X\}$. Consequently, there exists $\lambda^*_{\mu}\in
X$, such that the function $x\to J(x)-{{\mu}\over
{2}}\|x-\lambda^*_{\mu}\|^2$ has at least two global minima. 
By $(10)$ and the finite-dimensionality of $X$ again, the same function
satisfies the Palais-Smale condition, and so it admits at least three critical
points, thanks to Corolary 1 of [4]. Of course, this gives the conclusion, taking
$y_{\mu}=\lambda^*_{\mu}$.
\hfill $\bigtriangleup$\par
\medskip
REMARK 3.6. - Clearly, there are two situations in which Theorem 1.2 can immediately
be applied: when $\eta=+\infty$, and when $\eta>0$ and $\theta=0$. Note
that one has $\theta=0$ if, in particular, $J$ possesses at least two global
minima. \par
\medskip
REMARK 3.7. - It is also clear that under the same assumptions as those of Theorem 1.2
 but the finite-dimensionality of $X$, the conclusion is still true 
for every $\mu\in ]2\theta,2\eta[$ such that, for each $\lambda\in X$, the functional
$x\to J(x)-{{\mu}\over {2}}\|x-\lambda\|^2$ is weakly lower
semicontinuous and satisfies the Palais-Smale condition.\par
\medskip
{\it Proof of Theorem 1.3}. First of all, observe that $\varphi\in {\cal
G}$, with $y_0=0$. Indeed, let $V\subset Y$ be a neighbourhood of
$0$ and let $s>0$ be such that $B(0,s)\subseteq V$. Set
$$\alpha=\inf_{\|x\|=s}\varphi(x)\ .$$
Since dim$(Y)<\infty$, $\partial B(0,s)$ is compact and so $\alpha>0$.
Let $x\in Y$ with $\|x\|>s$. Let $S$ be the segment joining $0$ and $x$.
By convexity, we have
$$\varphi(z)\leq \varphi(x)$$
for all $z\in S$. Since $S$ meets $\partial B(0,s)$, we infer that
$\alpha\leq \varphi(x)$. Hence, we have
$$\alpha\leq \inf_{\|x\|>s}\varphi(x)\leq \inf_{x\in X\setminus V}\varphi(x)
\ .$$
Now, fix $x_0\in Y$ with $J'(x_0)\neq 0$. Taking into account that
$\varphi'(0)=0$, by continuity, we can choose $\sigma>0$ so that
$$\|\varphi'(\lambda)\|<\|J'(x)\|$$
 for all $(x,\lambda)\in B(x_0,\sigma)\times B(0,\sigma)$.
For each $(x,\lambda)\in Y\times Y$,
put
$$f(x,\lambda)=J(x)-\varphi(x-x_0-\lambda)\ .$$
Of course, we have
$$f'_x(x,\lambda)\neq 0$$
for all $(x,\lambda)\in B\left ( x_0,{{\sigma}\over {2}}\right ) 
\times B\left ( 0,{{\sigma}\over {2}}\right )$.
Next, since $J'$ is locally Lipschitzian at $x_0$  and $\varphi'$ is
locally Lipschitzian at $0$, there are $\rho\in \left 
] 0,{{\sigma}\over {2}}\right ] $ and
$L>0$ such that
$$\|f'_x(x,\lambda)-f'_x(y,\lambda)\|\leq L\|x-y\|$$
for all $x,y\in B(x_0,\rho)$, $\lambda\in B(0,\rho)$. Now,
fix $\lambda\in B(0,\rho)$. Denote by $\Gamma_{\lambda}$ the set
of all global minima of the restriction of the
function $x\to f(x,\lambda)+{{L}\over {2}}\|x-x_0\|^2$ to 
$B(x_0,\rho)$. Note that $x_0\not\in \Gamma_{\lambda}$ (since
$f'_x(x_0,\lambda)\neq 0$). As $f$ is
continuous, the multifunction $\lambda\to \Gamma_{\lambda}$ is
upper semicontinuous and so the function $\lambda\to \hbox
{\rm dist}(x_0,\Gamma_{\lambda})$ is lower semicontinuous. As
a consequence, by compactness, we have
$$\delta:=\inf_{\lambda\in B(0,\rho)}
\hbox {\rm dist}(x_0,\Gamma_{\lambda})>0\ .$$
At this point, from the proof of Theorem 1 of [6] it follows that, for
each $\lambda\in B(0,\rho)$ and each $r\in ]0,\delta[$, 
the restriction of function $f(\cdot,\lambda)$ to $B(x_0,r)$ has a
unique global minimum. Fix $r\in ]0,\delta[$. Apply Theorem 3.3
with $X=B(x_0,r)$, $\Lambda=Y$, ${\cal N}=\{B(x_0,r)\}$ and $\Psi(x,\lambda)=
x-x_0-\lambda$. With such choices, its conclusion does not hold with
$\mu=1$ (recall, in particular, that
$r<\rho$). This implies that $1\leq \theta(\varphi,\Psi,J)$ since the other
assumptions are satisfied. But the above inequality is just equivalent to
$$J(u_r)\leq J(x)-\varphi(x-u_r)$$
for all $x\in B(x_0,r)$, where $u_r$ is the unique global minimum of
$J_{|B(x_0,r)}$, and the proof is complete.\hfill $\bigtriangleup$\par
\medskip
Before proving Theorem 1.4, we point out another consequence of Theorem 3.1.\par
\medskip
THEOREM 3.4. - {\it Let $Y$ be a inner product space, and let
$I:X\to {\bf R}$, $\Phi:X\to Y$ and $\mu>0$ be such that
the function $x\to I(x)+\mu\|\Phi(x)\|^2$ has a global minimum.\par
Then, at least one of the following assertions holds:\par
\noindent
$(a)$\hskip 5pt for each filtering cover ${\cal N}$ of $X$, there exists
$A\in {\cal N}$ such that
$$\sup_{\lambda\in Y}\inf_{x\in A}(I(x)+\mu(2\langle\Phi(x),\lambda\rangle-
\|\lambda\|^2))<
\inf_{x\in A}\sup_{\lambda\in \Phi(A)}(I(x)+\mu(2\langle\Phi(x),\lambda\rangle-\|\lambda\|^2))\ ;$$
\noindent
$(b)$\hskip 5pt for each global minimum $u$ of $x\to I(x)+\mu\|\Phi(x)\|^2$,
one has
$$I(u)\leq I(x)+2\mu(\langle\Phi(x),\Phi(u)\rangle-\|\Phi(u)\|^2)$$
for all $x\in X$.}\par
\smallskip
PROOF. Take $\Lambda=Y$, $y_0=0$. For each $x\in X$, $y, \lambda\in Y$, set
$$\varphi(y)=\|y\|^2\ ,$$
$$\Psi(x,\lambda)=\Phi(x)-\lambda$$
and
$$J(x)=I(x)+\mu\|\Phi(x)\|^2\ .$$
So that
$$J(x)-\mu\varphi(\Psi(x,\lambda))=
I(x)+\mu(2\langle\Phi(x),\lambda\rangle -\|\lambda\|^2)\ .$$
With these choices, 
 $(b)$ is equivalent to the inequality
$$\mu\leq \theta(\varphi,\Psi,J)\ .$$
Now, the conclusion is a direct consequence of Theorem 3.1.
\hfill $\bigtriangleup$ \par
\medskip
A remarkable consequence of Theorem 3.4 is as follows:\par
\medskip
THEOREM 3.5. - {\it  Let $Y$ be a inner product space, and let
$I:X\to {\bf R}$, $\Phi:X\to Y$ be such that\par
\noindent
$(i)$\hskip 5pt $\sup_XI<+\infty$\ ;\par
\noindent
$(ii)$\hskip 5pt $\inf_{z\in \Phi(X)}\langle z,y\rangle=
-\infty$ for all $y\in \Phi(X)\setminus \{0\}$\ .\par
Under such hypotheses, for each $\mu>0$,
at least one of the following assertions holds:\par
\noindent
$(a)$\hskip 5pt for each filtering cover ${\cal N}$ of $X$, there exists
$A\in {\cal N}$ such that
$$\sup_{\lambda\in Y}\inf_{x\in A}(I(x)+\mu(2\langle\Phi(x),\lambda\rangle-
\|\lambda\|^2))<
\inf_{x\in A}\sup_{\lambda\in \Phi(A)}(I(x)+
\mu(2\langle\Phi(x),\lambda\rangle-\|\lambda\|^2))\ ;$$
$(b)$\hskip 5pt the set of all global minima in $X$
of the function $x\to I(x)+\mu\|\Phi(x)\|^2$ is contained in $\Phi^{-1}(0)$.
}\par
\smallskip
PROOF. Assume that $(a)$ does not hold. Then, we have to show that $(b)$
holds. So, let $u\in X$ be a global minimum
of the function $x\to I(x)+\mu\|\Phi(x)\|^2$. Then, by Theorem 3.4,  
 we have
$$I(u)\leq I(x)+2\mu(\langle\Phi(x),\Phi(u)\rangle-\|\Phi(u)\|^2) \eqno{(11)}$$
for all $x\in X$. Now, in view of $(i)$ and $(ii)$, if $\Phi(u)\neq 0$,
we would have
$$\inf_{x\in X}(I(x)+2\mu\langle\Phi(x),\Phi(u)\rangle)=-\infty$$
in contradiction with $(11)$. Therefore, $\Phi(u)=0$.
\hfill $\bigtriangleup$
\medskip
On the basis of Theorem 3.5, we can now give the\par
\medskip
{\it Proof of Theorem 1.4}. First of all, note that the function
$x\to I(x)+\langle\Phi(x),\lambda\rangle$ is convex in $X$ for all $\lambda\in Y$.
Now, set
$${\cal N}=\{\hbox {\rm conv}(B): B\subset X,\hskip 3pt
 B\hskip 3pt \hbox {\rm finite}\}\ .$$
Clearly, ${\cal N}$ is a filtering cover of $X$. Fix $A\in {\cal N}$.
Since $\Phi$ is affine, we have
$$\Phi(A)=\hbox {\rm conv}(\Phi(B))$$
where $B$ is a finite subset of $X$ such that conv$(B)=A$. So, $\Phi(A)$
is a convex and compact subset of $Y$. Taking into account that,
for each $x\in X$, the function $\lambda\to 
2\langle\Phi(x),\lambda\rangle
-\|\lambda\|^2$ is continuous and concave in $Y$, we can apply
Kneser's minimax theorem ([2]). So, thanks to it, we have
$$\sup_{\lambda\in \Phi(A)}\inf_{x\in A}
(I(x)+2\langle\Phi(x),\lambda\rangle -\|\lambda\|^2)=
\inf_{x\in A}\sup_{\lambda\in \Phi(A)}
(I(x)+2\langle\Phi(x),\lambda\rangle -\|\lambda\|^2)\ .$$
This shows that condition $(a)$ of Theorem 3.5 does not hold. So,
condition $(b)$ of the same theorem holds, and the conclusion is reached.
\hfill $\bigtriangleup$\par
\medskip
REMARK 3.8. - Let $(Y,\langle\cdot,\cdot\rangle)$ be a real inner product
space and let ${\cal R}$ denote the family of all subsets $S$ of $Y$,
containing more than one point, such that
$$\inf_{z\in S}\langle z,y\rangle=-\infty$$
for all $y\in S\setminus \{0\}$. In view of Theorems 3.5 and 1.4,
it would be interesting to characterize the sets
belonging to the family ${\cal R}$.
One of such characterizations, for convex sets, is suggested by Theorem 1.4 itself.
Actually, a convex set $S\subseteq Y$ belongs to ${\cal R}$ if and only
if, for each convex bounded above function $I:X\to {\bf R}$, the set
of all global minima in $S$ of the function $x\to I(x)+\|x\|^2$ is
either $\emptyset$ or $\{0\}$. Concerning the proof,
the "only if" part comes out directly from Theorem 1.4 taking $\Phi(x)=x$.
For the "if" part, it is enough to observe that, if we take $y\in
S\setminus \{0\}$, since
$y$ is the global minimum of
the function $x\to -2\langle x,y\rangle+\|x\|^2$, by assumption,
the convex function $x\to -2\langle x,y\rangle$ must be unbounded above in
$S$, and so $S\in {\cal R}$. In [13], J. Saint Raymond proved that,
when $Y$ is a Hilbert space, a closed convex subset of $Y$ (different
from $\{0\}$) belongs to the family ${\cal R}$ if and only if it is
a linear subspace. At the same time, he provided an example of a
closed convex set in a non-complete space $Y$ which belongs to ${\cal R}$
and does not contain $0$. So, it would be of particular interest to
find a geometric condition characterizing closed convex sets belonging
to ${\cal R}$ which, when $Y$ is complete, reduces to being a linear subspace.\par
\medskip
The proof of Theorem 1.5 is based on Theorem 3.5 again.\par
\medskip
{\it Proof of Theorem 1.5}. Let $X\subseteq L^p(T,E)$ be a decomposable set containing
the constant functions. Arguing by contradiction, assume that
$\hat u\in X$ is a global minimum of $J_{|X}$. For each $h\in L^p(T)$,
set
$$X_h=\{u\in X : \|u(t)\|\leq h(t)\hskip 3pt \hbox {\rm a.e.\hskip 3pt
in\hskip 3pt T}\}\ .$$
Clearly, $X_h$ is decomposable and
the family $\{X_h\}_{h\in L^p(T)}$ is a filtering covering of
$X$. For each $u\in X$, put
$$I(u)=\int_Tf(u(t))d\mu$$
and
$$\Phi(u)=\int_Tg(u(t))d\mu\ .$$
Note that our assumptions readily imply 
$$\sup_XI<+\infty,\hskip 3pt \inf_X\Phi=-\infty,\hskip 3pt
\sup_X\Phi=+\infty\ .\eqno{(12)}$$
Fix $h\in L^p(T)$. Thanks to the Lyapunov convexity theorem, $\Phi(X_h)$
is an interval, and so 
$$\Lambda_h:=\overline {\Phi(X_h)}$$
is a compact interval. Therefore, taking into account that the
function $$(u,\lambda)\to I(u)+2\lambda\Phi(u)-\lambda^2$$ is lower
semicontinuous in $L^p(T,E)$, and continuous and concave in ${\bf R}$,
we can apply Theorem 1.D of [7]. Thanks to it, we then have
$$\sup_{\lambda\in \Lambda_h}\inf_{u\in X_h}(I(u)+2\lambda\Phi(u)-\lambda^2)=
\inf_{u\in X_h}\sup_{\lambda\in \Lambda_h}(I(u)+2\lambda\Phi(u)-\lambda^2)\ .$$
Hence, condition $(a)$ of Theorem 3.5 does not hold, and so, in view
of $(12)$, condition $(b)$ of it must hold. 
Hence, $\Phi(\hat u)=0$,
and so $\hat u$ is a global minimum of $I$. This implies that $f$ has a
unique minimum, say $\xi_0$, and that $u(t)=\xi_0$ a.e. in $T$. Hence,
$g(\xi_0)=0$, against the assumptions.\hfill $\bigtriangleup$
\medskip
A further consequence of Theorem 3.1 is the following:
\par
\medskip
THEOREM 3.6. - {\it Let $I, \Phi:X\to {\bf R}$ and $\mu>0$ be such
that the function $I-\mu\Phi$ has a global minimum.\par
Then, at least one of the following assertions holds:\par
\noindent
$(a)$\hskip 5pt for each filtering cover ${\cal N}$ of $X$, there
exists $A\in {\cal N}$ such that
$$\sup_{\lambda\in {\bf R}}\inf_{x\in A}(I(x)-\mu(e^{\Phi(x)-\lambda}
+\lambda))<\inf_{x\in A}\sup_{\lambda\in \Phi(A)}
(I(x)-\mu(e^{\Phi(x)-\lambda}+\lambda))\ ;$$
\noindent
$(b)$\hskip 5pt for each global minimum $u$ of $I-\mu\Phi$, one has
$$I(u)\leq I(x)-\mu(e^{\Phi(x)-\Phi(u)}-1)$$
for all $x\in X$.}\par
\smallskip
PROOF. Take $\Lambda=Y={\bf R}$, $y_0=0$. For each
$x\in X$, $y, \lambda\in Y$, set
$$\varphi(y)=e^y-y-1\ ,$$
$$\Psi(x,\lambda)=\Phi(x)-\lambda$$
and
$$J(x)=I(x)-\mu \Phi(x)\ .$$
So that
$$J(x)-\mu\varphi(\Psi(x,\lambda))=I(x)-\mu(e^{\Phi(x)-\lambda}+\lambda - 1)\ .$$
With these choices, 
 $(b)$ is equivalent to the inequality
$$\mu\leq \theta(\varphi,\Psi,J)\ .$$
Now, the conclusion is a direct consequence of Theorem 3.1.
\hfill $\bigtriangleup$ \par
\medskip
Finally, we can give the\par
\medskip
{\it Proof of Theorem 1.6}. The proof is very similar to that of Theorem 1.4.
So, let ${\cal N}$ be as in that proof. Fix $A\in {\cal N}$. Since $P$ is
concave and positive, $P(A)$ is a bounded subset of ${\bf R}$ with positive
infimum. For each $x\in X$, put
$$\Phi(x)=\log P(x)\ .$$
Now, fix a compact real interval $[a,b]$ such that $\Phi(A)\subseteq [a,b]$.
Clearly, the function $(x,\lambda)\to I(x)-\mu (e^{\Phi(x)-\lambda}+\lambda)$ is convex
in $X$, and concave and continuous in ${\bf R}$. So, by Kneser's minimax theorem
again, we have
$$\sup_{\lambda\in [a,b]}\inf_{x\in A}(I(x)-\mu(e^{-\lambda}P(x)+\lambda))=
\inf_{x\in A}\sup_{\lambda\in [a,b]}(I(x)-\mu(e^{-\lambda}P(x)+\lambda))\ .$$
This shows that condition $(a)$ of Theorem 3.6 does not hold. So,
condition $(b)$ of the same theorem holds, and the conclusion is reached.
\hfill $\bigtriangleup$
\par
\vfill\eject
\centerline {\bf References}\par
\bigskip
\bigskip
\noindent
[1]\hskip 5pt V. L. KLEE, {\it Convexity of Chebyshev sets}, Math. Ann.,
{\bf 142} (1960/1961), 292-304.\par
\smallskip
\noindent
[2]\hskip 5pt H. KNESER, {\it Sur un th\'eor\`eme fondamental de la th\'eorie des jeux},
C. R. Acad. Sci. Paris, {\bf 234} (1952), 2418-2420.\par 
\smallskip
\noindent
[3]\hskip 5pt S. J. N. MOSCONI,
{\it A differential characterisation of the minimax inequality}, J. Convex
Anal., {\bf 19} (2012), 185-199.\par
\smallskip
\noindent
[4]\hskip 5pt P. PUCCI and J. SERRIN, {\it A mountain pass
theorem}, J. Differential Equations, {\bf 60} (1985), 142-149.\par
\smallskip
\noindent
[5]\hskip 5pt B. RICCERI, {\it On a three critical points theorem},
 Arch. Math. (Basel), {\bf 75} (2000), 220-226.\par
\smallskip
\noindent
[6]\hskip 5pt B. RICCERI, {\it The problem of minimizing locally a $C^2$ functional
around non-critical points is well-posed}, Proc. Amer. Math. Soc., {\bf 135}
(2007), 2187-2191.\par
\smallskip
\noindent
[7]\hskip 5pt B. RICCERI, {\it Recent advances in minimax theory and
applications}, in ``Pareto Optimality, Game Theory and Equilibria'', 
A. Chinchuluun, P.M. Pardalos, A. Migdalas, L. Pitsoulis eds., 23-52,
Springer, 2008.\par
\smallskip
\noindent
[8]\hskip 5pt B. RICCERI, {\it Nonlinear eigenvalue problems},
in ``Handbook of nonconvex analysis and applications'', D. Y Gao and
D. Motreanu eds., 543-595, International Press, 2010.\par
\smallskip
\noindent
[9]\hskip 5pt B. RICCERI, {\it Multiplicity of global minima for
parametrized functions}, Rend. Lincei Mat. Appl., {\bf 21} (2010),
47-57.\par
\smallskip
\noindent
[10]\hskip 5pt B. RICCERI, {\it
A multiplicity result for nonlocal problems involving nonlinearities with bounded primitive}, 
Stud. Univ. Babes-Bolyai Math., {\bf 55} (2010), 107-114.\par
\smallskip
\noindent
[11]\hskip 5pt B. RICCERI, {\it Addendum to ``A multiplicity result for
nonlocal problems involving nonlinearities with bounded primitive"},
Stud. Univ. Babes-Bolyai Math., {\bf 56} (2011), 173-174.\par
\smallskip
\noindent
[12]\hskip 5pt J. SAINT RAYMOND, {\it On a minimax theorem}, Arch. Math.
(Basel), {\bf 74} (2000), 432-437.\par
\smallskip
\noindent
[13]\hskip 5pt J. SAINT RAYMOND, Personal communication.\par
\smallskip
\noindent
[14]\hskip 5pt H. TUY, {\it A new topological minimax theorem with application}, J. Global Optim.,
 {\bf 50} (2011), 371-378. \par
\bigskip
\bigskip
Department of Mathematics\par
University of Catania\par
Viale A. Doria 6\par
95125 Catania, Italy\par
{\it e-mail address:} ricceri@dmi.unict.it

\bye
We now start a series of consequences and applications of
Theorem 1.\par
\medskip
THEOREM 2. - {\it Let $(H,\langle\cdot,\cdot\rangle)$ be
a real Hilbert space, $X$ a convex set in a topological vector
space, $J:X\to {\bf R}$ a lower semicontinuous
convex function and $T:X\to H$ a continuous affine operator with
bounded range.\par
Under such hypotheses, if $u\in X$ is a global minimum of
the function $x\to J(x)+{{1}\over {2}}\|T(x)\|^2$, then one has
$$J(u)\leq J(x)+\langle T(u),T(x)\rangle -\|T(u)\|^2$$
for all $x\in X$.}\par
\smallskip
PROOF. We apply Theorem 1 taking $Y=\Lambda=H$ with the
weak topology,
$$\varphi(y)={{1}\over {2}}\|y\|^2$$
for all $y\in H$, and
$$\Psi(x,\lambda)=T(x)-\lambda$$
for all $(x,\lambda)\in X\times H$. So that
$\lambda_x=T(x)$ for all $x\in X$. Finally, take
$$J(x)=J(x)+{{1}\over {2}}\|T(x)\|^2$$
for all $x\in X$. Note that the function
$x\to J(x)-\varphi(\Psi(x,\lambda))$ is lower semicontinuous and
convex
for all $\lambda\in H$, while the function
 $\lambda\to J(x)-\varphi(\Psi(x,\lambda))$ is upper
semicontinuous and concave for all $x\in X$. Then, since
$X$ is weakly compact, the classical Fan-Sion minimax
theorem, we have
$$\sup_{\lambda\in H}\inf_{x\in X}(J(x)-\varphi(\Psi(x,\lambda))=
\inf_{x\in X}\sup_{\lambda\in H}(J(x)-\varphi(\Psi(x,\lambda))\ .$$
That is to say, $(c)$ does not hold. Therefore, if $u$ is a global minimum of
$J$, then $u$ is also a global minimum of $x\to J(x)-\varphi(\Psi(x,\lambda_u))=
J(x)+{{1}\over {2}}\|T(x)\|^2-{{1}\over {2}}\|T(x)-T(u)\|^2$.
So
$$J(u)+{{1}\over {2}}\|T(u)\|^2\leq J(x)+{{1}\over {2}}\|T(x)\|^2-{{1}\over
{2}}\|T(x)-T(u)\|^2$$

ote also the following result that comes from a joint application
of Theorem 3.1 and Theorem 14 of []:\par
\medskip
THEOREM 3.6.- {\it Let $\varphi\in {\cal G}$, $\Psi\in {\cal H}$  
 and $J\in {\cal M}$. 
Moreover, assume that $X$ is a Hausdorff topological space, that
$\Lambda={\bf R}$ and that
$\varphi(\Psi(x,\cdot))$ is convex for each $x\in X$. Finally, let $\mu>\theta(\varphi,\Psi,J)$ be
such that $J(\cdot)-\mu\varphi(\Psi(\cdot,\lambda))$ is
inf-compact for each $\lambda\in
\hbox {\rm conv}(\Lambda_0)$.
\par
Under such hypotheses, there exists
a non-empty open set $A\subseteq \hbox {\rm conv}(\Lambda_0)$ such
that, for each $\lambda\in A$,
 $J(\cdot)-\mu\varphi(\Psi(\cdot,\lambda))$ has at least two
local minima.}\par
\medskip

\medskip
THEOREM G. - {\it Let $X$ be a real Hilbert space and let $J:X\to {\bf R}$ be a $C^1$ functional with locally
Lipschitzian derivative.\par
Then, for each $x_0\in X$ with $J'(x_0)\neq 0$, there exists $\delta>0$
such that, for every $r\in ]0,\delta[$, one has
$$\inf_{B(x_0,r)}J=\inf_{S(x_0,r)}J$$
and 
the problems of minimizing $J$
over $S(x_0,r)$ and over $B(x_0,r)$ are well-posed.}\par
\medskip

\vfill\eject
\centerline {\bf References}\par
\bigskip
\bigskip
\noindent
[4]\hskip 5pt B. RICCERI, {\it A three critical points theorem revisited},
Nonlinear Anal., {\bf 70} (2009), 3084-3089.\par
\smallskip
\noindent
[7]\hskip 5pt G. CORDARO, {\it On a minimax problem of Ricceri},
J. Inequal. Appl., {\bf 6} (2001), 261-285.\par
\smallskip
\noindent
\smallskip
\noindent
[10]\hskip 5pt B. RICCERI, {\it Sublevel sets and global minima of coercive
functionals and local minima of their perturbations}, J. Nonlinear Convex
Anal., {\bf 5} (2004), 157-168.\par
\smallskip
\noindent
[11]\hskip 5pt E. ZEIDLER, {\it Nonlinear functional analysis and its
applications}, vol. III, Springer-Verlag, 1985.\par
\smallskip
\noindent
[12]\hskip 5pt P. PUCCI and J. SERRIN, {\it A mountain pass theorem},
J. Differential Equations, {\bf 60} (1985), 142-149.\par
\par

\bye